\documentstyle[12pt]{article}

\def\Z{Z\!\!\!Z}
\let\ov\overline

\newtheorem{theorem}{Theorem}[section]

\newtheorem{lemma}[theorem]{Lemma}

\def\rem{\refstepcounter{theorem}\paragraph{Remark \thetheorem}}
\def\rems{\refstepcounter{theorem}\paragraph{Remarks \thetheorem}}
\def\proof{\paragraph{Proof}}

\makeatletter
\def\l@section{\@dottedtocline{1}{0em}{1.2em}}
\makeatother

\begin{document}
\title{Topology of Conic Bundles - II}
\author{Nitin Nitsure}
\date{24 February 1996}
\maketitle

\centerline{Tata Institute of Fundamental Research,
Mumbai 400 005, India.}
\centerline{e-mail: nitsure@tifrvax.tifr.res.in}

\centerline{\sl 1991 Mathematics subject classification :
14F20, 13A20, 12G05}

\begin{abstract}
\noindent{For} conic bundles on a smooth variety (over a field of
characteristic $\ne 2$) which degenerate into
pairs of distinct lines over geometric points of a smooth
divisor, we prove a theorem
which relates the Brauer class of the non-degenerate conic on the
complement of the divisor to the covering class (Kummer class)
of the 2-sheeted cover of the divisor defined by the degenerate
conic, via the Gysin homomorphism in etale cohomology. This
theorem is the algebro-geometric analogue of a topological
result proved earlier.
\end{abstract}

\section{Introduction}

Let $X$ be a smooth scheme over a field $F$ of
characteristic $\ne 2$, and let $C\to X$ be a conic bundle on
$X$, whose discriminant defines a smooth divisor $Y\subset X$ with
multiplicity $\tau$ (which, if $Y$ is not irreducible, will
consist of a positive integer $\tau_i$
for each component $Y_i$ of $Y$).
Let $\beta\in H^2(X,\mu_2)$ be the Brauer
class of the restriction of $C$ to $X-Y$, which is a $P^1$
fibration which is etale locally trivial (all cohomologies are
with respect to etale topology unless otherwise indicated).
Suppose that over each geometric point of $Y$, the fiber of $C$
consists of two distinct projective lines meeting at a point.
Therefore, the relative Hilbert scheme of lines in $C|Y \to Y$ is
a two sheeted finite etale cover $\tilde{Y}\to Y$. Let $\alpha \in
H^1(Y,\mu_2)$ be its covering class (`Kummer class').
By smoothness, we have a Gysin homomorphism $H^2(X-Y, \mu_2)\to
H^1(Y,\mu_2)$.  We prove here the following
\begin{theorem}
Under the Gysin homomorphism $H^2(X-Y, \mu_2)\to H^1(Y,\mu_2)$,
the Brauer class $\beta$ of the $P^1$ fibration on $X-Y$ maps to
$\tau\alpha$, where $\tau$ is the vanishing multiplicity of the
discriminant and $\alpha$ is the cohomology class of the two
sheeted finite etale cover $\tilde{Y}\to Y$.
\end{theorem}

A topological version of this result was proved in [Ni] for
topological conic bundles on manifolds, which
is equivalent for complex algebraic conic bundles to the
above result. This is because there is a natural isomorphism between
etale cohomology with finite constant coefficients (in this
case, coefficients $\Z/(2)$) and the corresponding singular
cohomology, which commutes with the two Gysins.

The theorem is proved below in two steps. In section 2, we
reduce it to proving a purely algebraic lemma (Lemma 2.1 below,
which we call as the `main lemma') over a discrete valuation ring.
In section 3, we prove the main lemma.

\rems (1) The main lemma is more transparent than
its topological counterpart in [Ni].

(2) After proving the main lemma, enquiries with algebraist
colleagues revealed that a more general lemma has already been
proved by Colliot-Th\'el\`ene and Ojanguren (see proposition 1.3
in [C-O]). Our proof is more geometric but less general.

(3) When the total space of $C$ is nonsingular,
it is known that (see [H-N] Proposition 1.4 or [Ne] Theorem 2)
$\beta$ is zero only if $\alpha$ is zero. This
follows from theorem 1.1 by taking $\tau=1$, though theorem 1.1
makes a stronger statement even in this case.

\section{Reduction to Main Lemma}
For basic definitions about conic bundles,
see for example Newstead [Ne].

It is clearly enough to prove the theorem for each connected
component of $X$, so we will assume $X$ to be connected.
If $Y=\cup _iY_i$ are the connected components of $Y$, then
by replacing $(X,Y)$ by $(X-\cup_{j\ne i}Y_j,Y_i)$ we are
reduced to the case where $Y$ also is connected. Hence we can
assume that both $X$ and $Y$ are irreducible.

Let the conic bundle $C\to X$ be defined via a rank
3 vector bundle $E$ on $X$, together with quadratic form
$q$ on $E$ with values in some line bundle $L$ on $X$.
The quadratic form $q$ is of rank 3 on $X-Y$ because we have a
non-degenerate conic over $X-Y$, and $q$ has rank 2 on $Y$ as we
have pairs of distinct lines over geometric points of $Y$.

Let $U = Spec R$
be an affine open subscheme of $X$ which intersects $Y$, such
that $E$ and $L$ are trivial on $U$, and $Y\cap U$ is defined by a
principal ideal $(\pi)$ in $R$. By injectivity of
$H^1(Y, \mu_2) \to H^1(Y\cap U, \mu_2)$, it is enough to prove the
theorem for $(U,Y\cap U)$ in place of $(X,Y)$.
Hence we can assume that the conic bundle is defined by an
explicit quadratic form $x^2-ay^2-bz^2$ on the ring $R$, where
$a$ is a unit in $R$, while $b\in (\pi)$ vanishes over $Y$.

Let $\eta$ be the generic point of $Y$, and let $A$ be the
discrete valuation ring $\O_{X,\eta}$. Let $k$ be the function
field of $Y$. The morphism $Spec(k) \to Y$ induces an injective
homomorphism $H^1(Y,\mu_2)\to H^1(k,\mu_2)$.
Hence the theorem follows from the following
lemma.

\begin{lemma}{\rm (Main Lemma) }
Let $F$ be a field of characteristic $\ne 2$.
Let $A$ be a discrete valuation ring which is the local ring
at the generic point of a smooth divisor in a smooth
$F$-variety. Let $K$ be the quotient field of $A$, let $k$ be
the residue field, and let $\nu :A-\{ 0\} \to \Z$ be the
discrete valuation.
Let $x^2-ay^2-bz^2$ be a quadratic form on $A$, with $a$ a unit
in $A$, and $b\ne 0$. Let $(a,b)\in H^2(K,\mu_2)$ be the Brauer class
(Hilbert symbol) of the quadratic form on $K$. Let $\ov{a}\in
k-\{ 0\}$ be the residue class of $a$, and let $\chi(\ov{a})\in
H^1(k,\mu_2)$ be the class of the two sheeted etale cover
$k(\ov{a}^{1/2})/k$ (Kummer character).
Then under the Gysin homomorphism $H^2(K,\mu_2)\to H^1(k,\mu_2)$, we
have (in additive notation)
$$(a,b) \mapsto \nu(b) \chi(\ov{a})$$
\end{lemma}

\section{Proof of the Main Lemma}

In the course of the proof below, we use the following elementary
facts which can be found for example in the textbook of Milne
[M]. If $S$ is a field or more generally a Henselian local ring,
then $Br(S)\to H^2(B, {\underline{GL}}_1)$ is an isomorphism, and
provided
the characteristic of the residue field is $\ne 2$, the
homomorphism $H^2(S,\mu_2)\to H^2(B,{\underline{GL}}_1)$ is injective.
Moreover the etale cohomology
of $S$ with coefficients in a smooth representable sheaf (for
example $\mu_2$ or ${\underline{PGL}}_2$) is isomorphic by
restriction to the corresponding etale
cohomology of the residue field of $S$.

We need two more fact, which are contained in the following two
remarks.

\rem\label{rem3.1}
Let $A$ be a henselian local ring and $B/A$ be a 2-sheeted
finite etale cover. If the image of $\gamma\in H^2(A,\mu_2)$ is
zero under the composite
$$H^2(A,\mu_2) \to H^2(B,\mu_2) \to H^2(L,\mu_2)$$
where $L$ is the quotient field of $B$, then $\gamma$ lies in the
image of the canonical (connecting) set map
$H^1(A,{\underline{PGL}}_2) \to H^2(A,\mu_2)$ for the following
reason. Any generic section of a Brauer-Severi
variety on $Spec(B)$ extends to a global section
by the valuative criterion of properness and so the map
$H^2(B,\mu_2) \to H^2(L,\mu_2)$ is injective, hence $\gamma$
maps to zero under $H^2(A,\mu_2) \to H^2(B,\mu_2)$. Hence
$\gamma$ is represented by an element (factor set)
of the group cohomology $H^2(Gal(B/A), \mu_2)$. As $Gal(B/A)$ is
of order 2, the factor set $\gamma$ defines an Azumaya algebra of
rank 2, showing $\gamma$ comes from $H^1(R,{\underline{PGL}}_2)$.

\rem\label{rem3.2} Let $K$ be a field of characteristic $\ne 2$,
let $a,b\in K-\{ 0\}$ such
that $a$ is not a square in $K$, and let $L=K[t]/(t^2-a)$.
As the conic in $P^2_K$ defined by $x^2-ay^2-bz^2=0$ has an
$L$-rational point, its Brauer class $(a,b)$ is
an element of the group cohomology set $H^1(Gal(L/K),
PGL_2(L))$. Let $Gal(L/K)=\{ 1,\sigma\}$. A 1-cocycle for
$Gal(L/K)$ with coefficients $PGL_2(L)$ therefore consists of
an element $g\in PGL_2(L)$ such that $g\sigma(g)=I$ in $PGL_2$.
(The corresponding `crossed homomorphism' $Gal(L/K)\to PGL_2(L)$
is defined by $\sigma\mapsto g$.)
If $g'\in GL_2(L)$ is an arbitrary lift of $g$, then there must
exist some $c\ne 0$ in $L$ such that $g'\sigma(g')=cI$ in
$GL_2(L)$ (which implies $c\in K$). In particular it can be seen
by using stereographic
projection from $(\sqrt{a}, 1,0)\in P^2_K(L)$ that
$(a,b)$ is represented in $H^1(Gal(L/K),PGL_2(L))$ by the
1-cocycle defined by
$$h=\pmatrix{
  & b \cr
1 &   \cr
}$$

We now prove the main lemma.

\proof If $b$ is a unit in $A$ (that is,
$\nu(b)=0$), then the $P^1$
bundle given by $x^2-ay^2-bz^2$ is defined over all of
$Spec(A)$, and so by exactness of the Gysin sequence
$$H^2(A,\mu_2)\to H^2(K,\mu_2) \to H^1(k,\mu_2)$$
the image of $(a,b)\in H^2(K,\mu_2)$ in $H^1(k,\mu_2)$ is also zero.
So we now assume $\ov{b}=0$ (which is anyway the case which interests
us). If $\nu(b) = 2n$ is even then as $\chi(\ov{a})$ is
2-torsion we have $\nu(b)\chi(\ov{a}) =0$.
Making the change of variable $z'= \pi^nz$ over $K$, we see that
$(a,b)=(a,b')\in H^2(K,\mu_2)$ where $b'=b/\pi^{2n}$ . As
$\nu(b')=0$, the argument above now
completes the proof when $\nu(b)$ is even. When $\nu(b)=2n+1$ is
odd, the same change of variables enables us to reduce to the case where
$\nu(b)=1$. Hence we assume from now on that $\nu(b)=1$.

Note that passing to the completion of $A$ does
not affect the residue field and gives a commutative square of
the Gysins. Hence it is enough to show the conclusion of the
lemma assuming $A$ is complete. If $\ov{a}\in k$ is a
square then (by the Henselian property of a complete local ring) $a$
will be a square in $A$ and hence in $K$, and conversely (this
uses the discrete valuation) if $a\in A$ is a square in $K$ then it is
a square in $A$ and hence $\ov{a}$ is a square in $k$. The conic
defined in $P^2_K$ by $x^2-ay^2-bz^2$ will have a $K$ rational
point $(a_1, 1, 0)$ whenever $a_1$ is a squareroot of $a$ in $K$,
showing its Brauer class over $K$ is zero. Hence if
$\chi(\ov{a})=0$ then $(a,b)=0$, and so the lemma holds. Hence
we now assume that $\ov{a}$ is not a square in $k$, and therefore
$a$ is not a square in $K$ or $A$. (As moreover $\nu(b)=1$, it
follows by remark 1.3 (3)  that
$(a,b)\ne 0$, though this will also follow from the argument below.)

Let $B=A[t]/(t-a^2)$ which is a 2-sheeted finite etale cover of $A$,
and let $L$ and $l$ be respectively the quotient field and the
residue field of $B$. Consider the commutative diagram where the
horizantal maps are the two Gysins.
$$\matrix{
H^2(K,\mu_2) & \to & H^1(k,\mu_2) \cr
\downarrow &       & \downarrow   \cr
H^2(L,\mu_2) & \to & H^1(l,\mu_2) \cr
}$$
The kernel of $H^1(k,\mu_2)\to H^1(l,\mu_2)$
is of order two, generated by $\chi(\ov{a})$. As $a$ has a
squareroot in $L$, $(a,b)\in Br(K)$ maps to zero in $Br(L)$.
Hence to show that $(a,b)\mapsto \chi(\ov{a})$, it is enough to
show that the image of $(a,b)$ in $H^1(k,\mu_2)$ is nonzero. By
exactness of the Gysin sequence, this will follow if we show
that $(a,b)$ does not lie in the image of $H^2(A,\mu_2)\to
H^2(K,\mu_2)$.

Suppose $(a,b)$ was the image of $\gamma \in H^2(A,\mu_2)$.
As $(a,b)$ restricts to zero over $L$, it follows from
remark \ref{rem3.1} that $\gamma$ lies in the image of some
$\gamma' \in H^1(A,{\underline{PGL}}_2)$ under
$H^1(A,{\underline{PGL}}_2)\to H^2(A,\mu_2)$.
As $\gamma'$ splits over $B$, it can be regarded as an element
of the group cohomology set $H^1(Gal(B/A),PGL_2(B))$, and hence
as in remark 3.2,
its 1-cocycle is represented by an element $g\in GL_2(B)$ such that
$$g\sigma(g) = eI$$
where $Gal(B/A)=\{ 1,\sigma\}$ and $e$ is a unit in $A$.

On the other hand, the group cohomology class of
$(a,b)\in H^1(Gal(L/K),PGL_2(L))$ can be represented by the matrix
$h$ given by remark 3.2.
If the cohomology class of $(g)$ were to map to the cohomology
class of $(h)$ in $H^1(Gal(L/K),PGL_2(L))$, then by definition
of group cohomology there would exist elements
$M\in GL_2(L)$ and $0\ne c\in L$ such that
$$h = cMg\sigma(M^{-1})$$
Applying $\sigma$, this would give
$$h=\sigma(h) = \sigma(c)\sigma(M)\sigma(g)M^{-1}$$
Multiplying the two equations and using $h^2=bI$ and
$g\sigma(g)=eI$, we would get
$$b/e = Norm(c)$$
But this is impossible as on one hand $e\in A$ is a unit while
$b\in A$ has valuation $\nu(b)=1$ so $\nu(b/e)=1$, and on the
other hand $\nu$ takes only even values on norms.

Hence the image of $(a,b)$ in $H^1(k,\mu_2)$ is nonzero.
This completes the proof of the main lemma and hence that of the
theorem.

\section*{References}
[C-O] Colliot-Th\'el\`ene J. L. and Ojanguren M. :
`Vari\'et\'es unirationnelles non rationnelles: au-del\`a de l'example
d'Artin et Mumford', Inventionnes Math. 97 (1989) 141-158.

[H-N] Hirschowitz, A. and Narasimhan, M. S. : `Fibr\'es de 't
Hooft speciaux et applications', in {\sl Enumerative geometry
and classical algebraic geometry}, Progress in Mathematics 24,
Birkhauser, 1982.

[M] Milne, J. S. : {\sl \'Etale Cohomology}, Princeton Univ.
Press, 1980.

[Ne] Newstead, P. E. : `Comparision theorems for conic bundles',
Math. Proc. Cambridge Philos. Soc. 90 (1981) 21-31.

[Ni] Nitsure, N. : `Topology of conic bundles', J. London Math. Soc.
(2) 35 (1987) 18-28.

\medskip

School of Mathematics, Tata Institute of Fundamental Research,
Homo Bhabha Road, Mumbai 400 005.
e-mail: nitsure@tifrvax.tifr.res.in

\end{document}